\documentclass[amsmath,amssymb,12pt,eps]{article}
\usepackage{color}

\newcommand{\beq}{\begin{equation}}
\newcommand{\eeq}{\end{equation}}

\newcommand{\beqs}{\begin{eqnarray}}
\newcommand{\eeqs}{\end{eqnarray}}

\usepackage{authblk}

\begin{document}

\title{\textsc{Generating Axially Symmetric Minimal Hyper-surfaces in ${\mathbf R}^{1,3}$}}

\author[1,2]{Jens Hoppe}
\author[3]{Jaigyoung Choe}
\author[4]{O. Teoman Turgut}
\affil[1]{Department of Mathematics, TUB
38106 Braunschweig,
Germany}
\affil[2]{Max Planck Institute for Mathematics, Vivatsgasse 7, 53111, Bonn, Germany}
\affil[3]{School of Mathematics, KIAS, 02455, Seoul, Republic of  Korea}
\affil[4]{Department of Physics, Bo\u{g}azi\c{c}i University, Bebek, 34342, \.{I}stanbul, T\"urkiye}

\maketitle

\abstract{It is shown that, somewhat similar to the case of classical
B\"acklund transformations for surfaces of constant negative curvature,
infinitely many axially symmetric minimal hypersurfaces in 4-dimensional
Minkowski-space can be obtained, in a non-trivial way, from any given one
by combining the scaling symmetries of the equations in light cone coordinates
with a non-obvious symmetry (the analogue of Bianchi’s original
transformation) - which can be shown to be involutive/correspond
to a space-reflection.}

\maketitle

\vskip 0.8 cm

Recently \cite{Hoppe1,Hoppe2} various signs of integrability for axially symmetric membranes in $4$-dimensional space time,
\begin{equation}
x^\nu(t,\varphi,\theta)=\pmatrix{t\cr r(t,\varphi)\matrix{\cos\theta\cr \sin \theta}\cr z(t,\varphi)}=\pmatrix{ \tau+\zeta/2\cr R(\tau,\mu)\matrix{\cos \Psi\cr \sin\Psi}\cr\tau-\zeta/2}=\tilde x^\nu(\tau, \mu,\Psi),
\end{equation}
were revealed, including the use of
\begin{equation}
{v'\over E}=\dot z,\ \ \ \ \dot v=r^2{z'\over E}
,\end{equation}
($\cdot  $ and $ '$ indicating derivatives with respect to $t$ and $\varphi$;
the constant $E$ sometimes put equal to $1$), where $r(t,\varphi)$ and $z(t,\varphi)$, describing the shape of the hypersurface, satisfy
\begin{eqnarray}
&\ &\dot r r'+\dot z z'=0, \ \ \ {\rm and} \ \ \ \dot r^2+\dot z^2 +r^2\Big( {r'^2+z'^2\over E^2}\Big)=1;\\
&\ & E^2 \ddot z= (r^2z')', \ \ \ \ \ \ \ \ \ \ \ E^2\ddot r=(r^2r')'-r(r'^2+z'^2) \nonumber
\end{eqnarray}
(the second line being implied by the first, provided $\dot r z'\neq r' \dot z$),
while in light-cone coordinates,
the first order equations of (3)
take the form
\begin{equation}
\dot \zeta={1\over 2} \Big(\dot R^2+{R^2 R'^2\over \eta^2}\Big), \ \ \ \ \zeta'=\dot R R'
,\end{equation}
implying
\begin{eqnarray}
&\ &\eta^2 \ddot R =R (R R')'  \\
&\ & \eta^2\ddot \zeta=(R^2\zeta')'
,\end{eqnarray}
with $R(\tau,\mu)=r(t,\varphi)$ and $\zeta=t-z$, $\cdot$ and $'$  denoting derivatives with respect to
\begin{equation}
\tau={1\over 2} (t+z(t,\varphi))\ \ \ \  {\rm and} \ \ \  \mu={1\over 2\eta} (E\varphi+v(t,\varphi)),
\end{equation}
 (2) implying
\begin{equation}
\eta\dot \mu={r^2\over E} \tau',   \ \ \ \ \eta {\mu'\over E}=\dot \tau,
\end{equation}
where $\eta$ is a constant (related to boosts in the $z$-direction, i. e. Minkowski rotations in the $(t,z)$-plane),  often put equal to $1$ just like $E$. Due to (5), respectively (6), (4) also implies the existence of a function $\kappa(\tau,\mu)$ satisfying
\begin{equation}
{\kappa'\over \eta}=\dot \zeta, \ \ \ \dot \kappa=R^2{\zeta'\over \eta}
.\end{equation}
(1), (8), and (9) are all generalizations of the Minkowski version of Cauchy-Riemann equations, i. e. for $r^2=1$ and $R^2=1$ all $6$ functions would satisfy the linear wave equation with respect to their independent variables ($t,\varphi$ for (1) and (8), $\tau,\mu$ for (9)).

Just as $g(f(w))$ is a holomorphic function of $w$, if $f$ and $g$ are, the following 'non-linear' (respectively $r$-dependent) generalization holds:\\
{\bf Theorem 1}:\\
Let $r(t,\varphi)$ be given, as well as $(x(t,\varphi), y(t,\varphi))$ and $(\phi(t,\varphi), T(t,\varphi))$ satisfy
\begin{eqnarray}
 &\ & x'=\dot y, \ \ \ \dot x=r^2 y'\\
&\ & \phi'=\dot T,\ \  \  \dot \phi=r^2 T'
\end{eqnarray}
$X(T, \phi)=x(t,\varphi)$ and $Y(T,\phi)=y(t,\varphi)$
(with $R(T,\phi)=r(t,\varphi)$)
 will then satisfy
\begin{equation}
X'=\dot Y, \ \ \ \dot X=R^2 Y'
\end{equation}
provided the change of variables $(t,\varphi)\mapsto (T,\phi)$ is invertible i. e. $ \dot T^2- r^2 T'^2\neq 0$.\\
{\it Proof:} We note that,
\begin{equation}
\partial_t=\dot T \partial_T+\dot \phi\partial_\phi, \ \ \ \partial_\varphi=T'\partial_T+\phi'\partial_\phi\, .
\end{equation}
So, using (11), (10) implies
$$
\pmatrix{\dot T & T'\cr r^2T' & \dot T} \pmatrix{ X'-\dot Y\cr \dot X-R^2Y'}=\pmatrix{ 0\cr 0},
$$
hence (12) follows.
For the construction of infinitely many non-trivial solutions of (5), the following observation will be crucial:\\
{\bf Theorem 2}:\\
Given (5), respectively (4) and (9),
$\rho(\zeta, \kappa)=R(\tau,\mu)$ will satisfy $\ddot \rho=\rho (\rho \rho')'$ i.e. a new solution of (5) can be generated by rewriting a starting solution $R$ in terms of $\zeta,\kappa$
(obtained by solving (4) and  (9)).\\
{\it Proof:} Note that
\begin{equation}
\pmatrix{\tau_\zeta& \tau_\kappa\cr \mu_\zeta& \mu_\kappa} = \pmatrix{ \dot \zeta& \zeta'\cr \dot \kappa & \kappa'}^{-1} ={1\over \delta} \pmatrix{ \kappa'& -\zeta'\cr
-\dot \kappa& \dot \zeta},
\end{equation}
where $\delta=\eta (\dot \zeta^2-R^2 {\zeta'^2\over \eta^2})={\eta\over 4} ( \dot R^2-R^2 {R'^2\over \eta^2})^2=:\eta {\mathcal L}^2$,
gives
\begin{eqnarray}
\dot \rho&=&\rho_\zeta= (\tau_\zeta \partial_\tau+\mu_\zeta \partial_\mu)R=...={\dot R\over {\mathcal L}}\nonumber\\
\rho'&=& \rho_\kappa= ( \tau_\kappa\partial_\tau+ \mu_\kappa\partial_\mu)R=...=-{R'\over \eta {\mathcal L}}
.\end{eqnarray}
A lengthy, but  straigthforward, calculation of $\delta\ddot \rho=(\kappa'\partial_\tau-\dot \kappa\partial_\mu) {\dot R\over {\mathcal L}}$ and
$\delta \rho''=-(-\zeta'\partial_\tau+\dot \zeta\partial_\mu) {R'\over \eta {\mathcal L}}$ then gives the desired result. Note also that (15) implies
\begin{equation}
{1\over 2} (\dot \rho^2-\rho^2 \rho'^2)\cdot {1\over 2} (\dot R^2 - {R^2 R'^2 \over \eta^2})=1.
\end{equation}
 It would at first sight be tempting to apply Theorem 2 directly multiple times, to obtain infinitely many solutions from a given one. However,  (16),  as well as
$(\zeta,\tau)$, as functions of $t, \varphi$, related to $(\tau, \mu)$ by the simple reflection $z\to -z$ (the factors of $2$ being irrelevant, respectively cancelling) indicate that the transformation is involutive. Indeed, considering
\begin{equation}
X'=\dot Y={1\over 2} (\dot \rho^2 +\rho^2\rho'^2) \ \ \ {\rm and} \ \ \ \dot X=\rho^2 Y'=\rho^2\dot \rho \rho'
\end{equation}
in $\tau, \mu$ coordinates, using ($\eta=1$),
\begin{equation}
\partial_\zeta={1\over {\mathcal L}^2} (\kappa'\partial_\tau-\dot \kappa \partial_\mu), \ \ \
\partial_\kappa={1\over {\mathcal L }^2} ( -\zeta'\partial_\tau+\dot \zeta \partial_\mu)
,\end{equation}
gives
$$
\kappa'\dot x-\dot\kappa x'=R^2(-\zeta'\dot y+\dot \zeta y') ={\mathcal L}^2 \rho^2\dot \rho\rho'=-R^2\dot R R'
$$
$$
-\zeta'\dot x +\dot \zeta x'=\kappa' \dot y -\dot \kappa y'={1\over 2} {\mathcal L}^2(\dot \rho^2 +\rho^2 \rho'^2) ={1\over 2} (\dot R^2 +R^2 R'^2)
$$
which implies
$$
\pmatrix{ \dot x \cr x'}=\pmatrix{0\cr 1}, \ \ \ \ \ \pmatrix{ \dot y\cr y'}=\pmatrix{1\cr 0}
,$$
i. e. up to some trivial additive constants
\begin{equation}
x(\tau,\mu)=\mu, \ \ \ y(\tau,\mu)=\tau
.\end{equation}
Just as with the original B\"acklund transformation \cite{back-4} which (as observed by Lie) were a combination of  Bianchi's original transformation \cite{bianchi-5} and trivial scaling symmetries ( see \cite{transform-6} for a nice discusion of $K=-1$ surfaces) one may combine Theorem 2 with the ('trivial') observation that
\begin{equation}
R_{\alpha,\gamma}(\tau,\mu)=\alpha R(\alpha\gamma\tau, \gamma \mu)
\end{equation}
satisfies (5), if $R$ does, to obtain infinitely many solutions of (5) which are \underline{not} simply results of (20)
(namely applying theorem 2 \underline{and} (20) infinitely many times, always in pairs).

As nice as this seems, explicit non-trivial elementary examples of the construction are not easy to come by:
\begin{equation}
R(\tau,\mu)=\sqrt{2}{\mu\over \tau},\ \ \ \ \zeta=-{\mu^2\over \tau^3}=-{R^2\over 2\tau}, \ {\rm and}\ \ \kappa={\mu^3\over \tau^4}
\end{equation}
is invariant under both (20) and the involutary transformation of theorem 2,
\begin{equation}
R\mapsto R^*(\tau,\mu):=\rho(\tau,\mu)
\end{equation}
as ${\kappa\over \zeta}=-{\mu\over \tau}$, hence $\rho(\zeta,\tau)=-\sqrt{2} {\kappa\over \zeta}$, i.e. $R^*=-R$ (which of course  is always a solution, to be identified with $R$  as the radius of the rotated curve should, by definition, be positive.) (21) corresponds to ${\mathcal M}_3=\{ x^\nu \in {\mathbf R}^{1,3} |
t^2+x^2+y^2=z^2\}$, moving hyperboloids, which in orthonormal parametrization were shown in \cite{Ref7} to read
\begin{equation}
z(t,\varphi)=\pm t \sqrt{ {\sqrt{ 1+ 8{\varphi^2\over t^4}} }+1\over 2}, \ \ \ \  r(t,\varphi)=|t| \sqrt{ {\sqrt{ 1+ 8{\varphi^2\over t^4}}}+1\over 2},
 \end{equation}
solving (2), hence also the second-order equations ($E=1$)
\begin{equation}
E^2\ddot r=(r^2r')'-r(r'^2+z'^2), \ \ E^2\ddot z=(r^2z')'.
\end{equation}
On the other hand, applying Theorem 2 to the minimal hypersurface
\begin{equation}
R=\sqrt{2} {\sqrt{\mu^2+\epsilon}\over \tau}, \ \ \ \zeta=-{\mu^2+\epsilon/3\over \tau^3}, \ \ \ \kappa={\mu^3+\epsilon\mu\over \tau^4}
,\end{equation}
(where $\epsilon$  is a  constant), whose level set form,
\begin{equation}
{\mathcal M}_3=\{ x^\nu\in {\mathbf R}^{1,3} | (t^2+x^2+y^2-z^2)(t+z)^2={16\epsilon\over 3}\}
\end{equation}
was, 30 years after Dirac's spherically symmetric solution \cite{dirac} the first non-trivial, polynomial solution \cite{Ref9}, is already very difficult to fully work out in detail, as one has to solve polynomial equations of high degree(s) to obtain $\mu(\zeta,\kappa)$ and $\tau(\zeta, \kappa)$ and hence $\rho(\zeta,\kappa)$, from (25) --or to find an explicit orthonormal parametrization of (26) by inverting
\begin{equation}
t=\tau-{\mu^2+\epsilon/3\over 2\tau^3}\ \ \ \ E\varphi=\mu+{\mu^3+\epsilon\mu\over 2\tau^4}\, .
\end{equation}
Before going to the next example,
note that (20) will give
\begin{equation}
\zeta_{\alpha,\gamma}(\tau,\mu)=\alpha^3\gamma\, \zeta(\alpha\gamma\tau, \gamma\mu),\ \ \kappa_{\alpha,\gamma}(\tau,\mu)=\alpha^4\gamma\,\kappa(\alpha\gamma\tau, \gamma\mu)
\end{equation}
and the effect of scale transformations in $(t,\varphi)$ parametrization can be calculated according to
\begin{eqnarray}
t&=&\tau+{\alpha^3\gamma\over 2}\, \zeta(\alpha\gamma\tau, \gamma\mu)\nonumber\\
E\varphi&=& \mu+ {\alpha^4\gamma\over 2}\,  \kappa(\alpha\gamma\tau, \gamma\mu)
\end{eqnarray}
the inversion of which gives $\tau=\tau_{\alpha,\gamma} (t,\varphi), \mu=\mu_{\alpha,\gamma}(t,\varphi)$ from which
\begin{equation}
z_{\alpha, \gamma}(t,\varphi)=\tau-{1\over 2} \zeta_{\alpha,\gamma},\ \ \ {\rm and} \ \ r_{\alpha,\gamma}(t,\varphi)=R_{\alpha,\gamma} (\tau,\mu)
\end{equation}
can in principle be calculated.

Consider now, as a third  (class of) example(s)
\begin{equation}
R=\tau\sqrt{\mu},\ \ \ \zeta={\mu \tau\over 2}+{\tau^5\over 40},\ \ \ \kappa={\tau^4\over 8}\mu+{\mu^2\over 4}
\end{equation}
respectively their scaled versions,
\begin{eqnarray}
R_{\alpha,\gamma}&=&R_{\beta=\alpha^2\gamma^{3/2}}(\tau,\mu)=\beta\tau\sqrt{\mu}\nonumber\\
\zeta_\beta&=&\beta^2({\mu\tau\over 2}+\beta^2{\tau^5\over 40})={R_\beta^2\over 2\tau}+{\beta^4\tau^5\over 40}\nonumber\\
\kappa_\beta&=& \beta^4{\tau^4\mu\over 8}+\beta^2{\mu^2\over 4}
\end{eqnarray}
corresponding to
\begin{equation}
{\mathcal M}=\{ x^\nu\in {\mathbf R}^{1,3} |\  t^2-x^2-y^2-z^2=C(t+z)^6\},\ \ {\rm with} \ \ C={\beta^4\over 1280}>0
\end{equation}
(discussed in \cite{Hoppe2,Hoppe1} and `quantized' in \cite{Ref10}). In order to obtain $R^*$, respectively $\rho$, by writing $R$ in (31)/(32)
as a function of $\zeta$ and $\kappa$ one would have to solve polynomial equations like ($\beta=1$)
\begin{equation}
(2\zeta-{\tau^5\over 20}) ( {9\over 20 } \tau^5+2\zeta)=4\tau^2\kappa.
\end{equation}
As a final example, consider
\begin{equation}
R=\sqrt{2}F(\tau)\mu, \ \ \ \zeta=F\dot F \mu^2,\ \ \ \kappa=F^4\mu^3+{1\over 3} \mu^3
\end{equation}
where $F$ is an elliptic function satisfying
\begin{equation}
\dot F^2=F^4+1\, .
\end{equation}
Here one can compute  that
\begin{equation}
\rho(\zeta,\kappa)=\zeta^{1/2} F(\kappa \zeta^{-3/2})
\end{equation}
showing that the involutary transformation $R\mapsto \rho$ of Theorem 2 maps two self-similar solutions (with different exponents) onto each other.

In order to appreciate the non-triviality of implementing the $z\mapsto -z$ symmetry on solutions $R(\tau, \mu), \zeta(\tau, \mu)$ of (4),(5) and (6),  consider for example (25)/(26): although Theorem-2 is difficult to apply explicitly (as for that $\mu$ and $\tau$ would be needed as explicit functions of $\zeta$ and $\kappa$), one may observe that $\zeta^4/\kappa^3$ depends only on $\mu$ hence $\mu=\mu(\zeta \kappa^{-3/4})$, implying $\tau=\zeta^{-1/3} h(\zeta \kappa^{-3/4})$ and $\rho=\zeta^{1/3} f(\zeta \kappa^{-3/4})$. Theorem 2 thus predicts a new solution $\tilde R(\tau,\mu)$, that is of the form
\begin{equation}
\tilde R(\tau,\mu)=\tau^{1/3} f(\mu\tau^{-3/4}:=\xi)
.\end{equation}
Both (4), and the $z\mapsto -z$ variant of (26),
\begin{equation}
(2\tau\tilde \zeta+\tilde R^2)\tilde \zeta^2={16\epsilon\over 3}=C\ (<0);
\end{equation}
then imply/suggest that
\begin{equation}
\tilde\zeta(\tau, \mu)=\tau^{-1/3} g(\xi)
.\end{equation}
Indeed, (38)/(40) consistently reduces  (39) and (4) to
\begin{equation}
2g^3+f^2 g^2= C\, (<0)
\end{equation}
and the 2 ODE's
\begin{eqnarray}
g'&=& f' ({1\over 3} f+\xi f')\\
2(\xi g'-{g\over 3})&=& ({1\over 3} f+\xi f')^2+ {9\over 16} \xi^{14/3} f^2 (f')^2.
\end{eqnarray}
Using (from (41))
\begin{equation}
f=-{1\over g} \sqrt{ C-2g^3}, \ \ \ \ f'={g'\over g^2\sqrt{C-2g^3}}(C+g^3)
,\end{equation}
(42) can be written as
\begin{eqnarray}
-{d\xi\over \xi} &=&-3{dg\over g} {(C+g^3)^2\over (C+4g^3)(C-2g^2)}\stackrel{(z=-g^3)}{=} {dz\over z} {(z-C)^2\over (2z+C)(4z-C)}\nonumber\\
&=& dz \Big(-{1\over z} + {3/4\over (z+ C/2)}+ {3/8\over (z-C/4)}\Big)
\end{eqnarray}
implying
\begin{equation}
 E^2 \xi= {z\over (z+C/2)^{3/4} (z-C/4)^{3/8}}.
\end{equation}
(43), on the other hand gives
\begin{eqnarray}
&\ &6(C+4g^3) (C-2g^3) g^3 -6g^3 (C+g^3)^2-(C-2g^3) (C+g^3)^2-(C+4g^3)(C-2g)\nonumber\\
&\ & \ \ \ \ \ \ \ \ \ \ \ \ \ \ \ \ \ \ \ \   +2(C+g^3)(C+4g^3)(C-2g^3)=
  {9\over 16}{\xi^{8/3}\over g^2} (C-2g^3)^2(C+4g^3)
\end{eqnarray}
with a (at first sight `weird') factor ${\xi^{8/3}\over g^2}$, which however,  using (46) (with $E^2=1$) is exactly what is needed to make (43) consistent with (42)/(41) (both sides of (47) are equal to $-9z^2(C-4z)$,  provided the integration constant $E^2$ is chosen to be equal to 1).

So, while (25) can not be explicitly inverted, the solution that is obtained from it via Theorem 2, respectively $z\mapsto -z$ in (26), can be obtained `almost explicitly', the function $g(\xi)$ in (40) solving a polynomial equation of degree 27,
\begin{equation}
\xi^8 (g^3+C/4)^3(g^3-C/2)^6=g^{24}
.\end{equation}
\vskip 0.5 cm

\noindent {\Large {\bf Acknowledgement}}\\
J.C is supported in part by Korea NRF-2018R1A2B6004262.

\end{document}